\documentstyle[11pt]{article}
\begin{document}
\title{{\bf Extremal subsets of $\{1,...,n\}$ avoiding solutions to 
linear equations in three variables}}
\author{{\sc Peter Hegarty} \\ Chalmers University of Technology and 
Gothenburg University, Sweden \\ hegarty@math.chalmers.se}
\date{June 19, 2007}
\maketitle

\begin{abstract}

We refine previous results to provide examples, and in some cases
precise classifications, of extremal subsets of $\{1,...,n\}$ containing 
no solutions to a wide class of non-invariant, homogeneous linear equations
in three variables, i.e.: equations of the form $ax+by=cz$ with 
$a+b \neq c$.

\end{abstract}

\section{Introduction}

A well-known problem in combinatorial number theory is that of locating
extremal subsets of $\{1,...,n\}$ which contain no non-trivial solutions to 
a given linear equation 
\begin{equation}
{\cal L} : \;\;\; a_{1}x_{1} + \cdots + a_{n} x_{n} = b,
\end{equation}
where $a_{1},...,a_{n},b \in {\hbox{{\bf Z}}}$ and their GCD is one. 
Most of the best-known
work concerns just three individual, homogeneous equations 
\begin{eqnarray*}
{\cal L}_{1} : \;\;\; x_{1} + x_{2} = 2x_{3}, \\
{\cal L}_{2} : \;\;\; x_{1}+x_{2} = x_{3} + x_{4}, \\
{\cal L}_{3} : \;\;\; x_{1}+x_{2} = x_{3},
\end{eqnarray*}
where the corresponding subsets are referred to, respectively, as 
sets without arithmetic progressions, Sidon sets and sum-free sets. The idea
to consider arbitrary linear equations ${\cal L}$ was first 
enunciated explicitly in a pair of articles by Ruzsa in the mid-1990s 
[10] [11]. The only earlier reference of note would appear to be a 
paper of Lucht [6] concerning homogeneous equations in three 
variables, though Lucht's article was only concerned with subsets of 
{\bf N}. Following Ruzsa, denote by $r_{{\cal L}}(n)$ the maximum size of
a subset of $\{1,...,n\}$ which contains no non-trivial solutions
to a given equation $\cal L$ (the precise meaning of the words 
$\lq$non-trivial' is given in [10]). When considering arbitrary $\cal L$, one
begins by observing a basic distinction between those $\cal L$ which 
are translation-invariant $(\sum a_{i} = b = 0)$ and those which are not, 
namely : for the former it is always the case that $r_{{\cal L}}(n) = o(n)$, 
a fact which follows easily from Szemer\'{e}di's famous theorem, whereas 
for the latter $r_{{\cal L}}(n) = \Omega(n)$ always. 
\par This paper is concerned with non-invariant, homogeneous equations only.
For simplicity the words $\lq$(linear) equation' will, for the remainder of 
the article, 
be assumed to refer to those equations with these extra properties,
though some of our initial observations also apply in the inhomogeneous
setting. We shall also employ the concise formulation
$\lq$$A$ avoids $\cal L$' to indicate that a set $A$ of 
positive integers contains no solutions to the equation $\cal L$. Finally
we will employ the interval notation $[\alpha,\beta] := 
\{x \in {\hbox{{\bf Z}}} : \alpha \leq x \leq \beta\}$, and similarly for
open intervals. 
\\
\\ 
As Ruzsa observed, given an equation
${\cal L} : \sum a_{i}x_{i} = 0$, there are two basic ways to exhibit the fact 
that $r_{{\cal L}}(n) = \Omega(n)$ :
\\
\\
{\bf I}. Let $s := \sum a_{i}$ so $s \neq 0$. Let $q$ be any positive 
integer not dividing $s$ and let $A := \{x \in {\hbox{{\bf N}}} : x 
\equiv 1 \; ({\hbox{mod $q$}})\}$. Then $A$ avoids $\cal L$ and 
$|A \cap [1,n]| \geq \lfloor n/q \rfloor$. 
\\
\\
{\bf II}. Set 
\begin{eqnarray*}
s_{+} := \sum_{a_{i} > 0} a_{i}, \;\;\; s_{-} := \sum_{a_{i} < 0} |a_{i}|
\end{eqnarray*}
and assume without loss of generality that $s_{+} > s_{-}$. For a fixed
$n > 0$ let $A := \left( \frac{s_{-}}{s_{+}} n, n \right]$. Then $A$ 
avoids $\cal L$ and has size $\Omega(n)$. 
\\
\\
As in [11], set $\lambda_{0,{\cal L}} := \limsup_{n \rightarrow \infty} 
{r_{{\cal L}}(n) \over n}$. Ruzsa asked whether the above two constructions
were the prototypes for extremal $\cal L$-avoiding sets in the sense 
that $\lambda_{0} = \max \{ \rho, {s_{+} - s_{-} \over s_{+}} \}$, where the
quantity $\rho$ is defined as follows : for each $m > 0$ let $\rho_{m}\cdot m$
be the maximum size of a subset of $[1,m]$ which contains no solutions to
$\cal L$ modulo $m$. Then $\rho := \sup_{m} \rho_{m}$. 
\par As illustrated by Schoen [12], the answer to Ruzsa's question is no. 
But from what is currently known, it seems that for many equations
something not much more complicated holds. One observes that 
the construction {\bf II} above can be modified into something more
general :
\\
\\
{\bf II}$^{\prime}$. For a given ${\cal L} : \sum a_{i}x_{i} = 0$, 
let notation be as 
above and let $a$ denote the smallest absolute value of a negative
coefficient $a_{i}$. Now for fixed $k,n > 0$ and $\xi \in [1,n]$ set 
\begin{eqnarray*}
A_{n,k,\xi} := \bigsqcup_{j=1}^{k-1} \left( {s_{-} \over s_{+}} n_{j}, n_{j} \right]
\sqcup [\xi,n_{k}],
\end{eqnarray*}
where $n_{1},...,n_{k}$ is any sequence of integers satisfying the recurrence
\begin{equation}
n_{1} = n, \;\;\; {s_{-} \over s_{+}} n_{k} < \xi \leq n_{k}, \;\;\; 
s_{+}n_{j+1} \leq a{s_{-} \over s_{+}}n_{j} + (s_{-} - a)\xi, 
\;\; j = 1,...,k-1.
\end{equation}
Clearly for (2) to have any solution we will need to have $k = O(\log n)$. 
Assuming a solution exists, the set $A_{n,k,\xi}$ avoids $\cal L$ and 
\begin{eqnarray*}            
|A_{n,k,\xi}| = (1+o(1)) \cdot \left[ {s_{+} - s_{-} \over s_{+}} 
\sum_{j=1}^{k-1} n_{j} + (n_{k} - \xi + 1) \right].
\end{eqnarray*}
The important special case is when $\xi = 1 + \lfloor {s_{-} \over s_{+}}
n_{k} \rfloor$. Then Ruzsa's question can be replaced by the following :
\\
\\
{\bf Question} {\em Is it always the case that 
\begin{eqnarray*}
\lambda_{0} = \max \left\{ \rho, \sup_{n,k,\xi} {|A_{n,k,\xi}| \over n} \right\},
\end{eqnarray*}
where the supremum ranges over all triples $n,k,\xi$ for which (2) has 
a solution where $\xi = 1 + \lfloor {s_{-} \over s_{+}} n_{k} \rfloor$ ?}
\\
\\
We will give examples in Section 3 
which show that the answer to this question is 
still no : we are not aware of any in the existing literature. However,
existing results strongly suggest that the answer is very often yes : \\
see [1] [2] [3] [8] [9] 
for example, plus further results in this paper. Also, in our
counterexamples the extremal sets are a pretty obvious hybrid between the two 
alternatives which the question offers. We think that our question is thus 
a good foundation for further research in this area.
\\
\\
We now give 
a closer overwiew of the results in this paper. To identify extremal 
$\cal L$-avoiding sets and compute $\lambda_{0,{\cal L}}$ for arbitrary $\cal L$
seems a very daunting task, so an obvious strategy is to study equations in 
a fixed number of variables. One variable is utterly 
trivial and two only slightly less so. I have not been able to locate the
following statement anywhere in the literature, however (though 
see for example [5], pp.30-34), so include it for
completeness :
\\
\\
{\bf Proposition} {\em Consider the equation ${\cal L} : ax=by$ where $a > b$ 
and GCD$(a,b) = 1$. For every $n > 0$ an extremal $\cal L$-avoiding subset
of $[1,n]$ is obtained by running through the numbers from $1$ to
$n$ and choosing greedily. This yields the extremal subset 
\begin{eqnarray*}
A := \{u \cdot a^{2i} : i \geq 0 \; {\hbox{and}} \; a \dagger u \}
\end{eqnarray*}
of {\bf N}. For each $n > 0$ a complete description of the extremal 
$\cal L$-avoiding subsets of $[1,n]$ is given as follows :
\\
\\
{\sc Case I} : $b = 1$.
\\
\\
For each $u \in [1,n]$ such that $u \dagger a$, let $\alpha$ be the largest
integer such that $u \cdot a^{\alpha} \leq n$. Then an extremal
set contains exactly $\lceil \alpha/2 \rceil$ of the numbers 
$u \cdot a^{i}$, for $0 \leq i \leq \alpha$, and no two numbers 
$u \cdot a^{i}$ and $u \cdot a^{i+1}$. 
\\
\\
{\sc Case II} : $b > 1$. 
\\
\\
For each $u \in [1,n]$ divisible by neither $a$ nor $b$ and each non-negative
integer $\alpha$ such that $u \cdot a^{\alpha} \leq n$, an extremal set 
contains exactly $\lceil \alpha/2 \rceil$ of the numbers 
$u \cdot b^{i} \cdot a^{\alpha -i}$, for $0 \leq i \leq \alpha$, 
and no two numbers
$u \cdot b^{i} \cdot a^{\alpha - i}$ and $u \cdot b^{i+1} \cdot 
a^{\alpha - i -1}$.}
\\
\\
Note that the proposition implies in particular that $\lambda_{0} = \rho$ for
any equation in two variables. For three variables things get more interesting
and a number of papers have been entirely devoted to this situation, see  
[1] [2] [4] [6] [7] plus the multitude of papers on sum-free sets, of
which the most directly relevant is probably [3]. The combined results of 
[1], [2] and [3] give, in principle, a complete classification of the 
extremal $\cal L$-avoiding subsets of $[1,n]$, for every $n > 0$,  
and ${\cal L} : x+y=cz$ for any $c \neq 2$. Of particular interest for us are 
the results of [1]. There it is shown that for every $c \geq 4$ and $n >>_{c}
0$, a set $A_{n,3}$ of type {\bf II}$^{\prime}$ is extremal, namely 
\begin{eqnarray*}
A_{n,3} := \bigsqcup_{j=1}^{3} \left( {2 \over c} n_{j}, n_{j} \right],
\end{eqnarray*}
where 
\begin{eqnarray*}
n_{1} = n, \;\;\;\; n_{j+1} = \lfloor {(1+\lfloor \frac{2}{c} n_{j} \rfloor) + 
(1+\lfloor \frac{2}{c} n_{3} \rfloor) \over c} \rfloor, \;\; j = 1,2.
\end{eqnarray*}
Moreover it is shown that, for all $n >>_{c} 0$, there are only a bounded
number of extremal sets, all of whose symmetric differences with $A_{n,3}$
consist of a bounded number of elements (both bounds are independent of
both $n$ and $c$). 
\par These results were partly extended in [4]. Here the authors 
considered equations ${\cal L} : ax+by=cz$ in two families :
\\
\\
{\sc Family I} : $a=1 < b$.
\\
{\sc Family II} : $a=b$, GCD$(b,c) = 1$.
\\
\\
For {\sc Family I} equations their main result is that, when 
$c > 2(b+1)^{2} - {\hbox{GCD}}(b+1,c)$ then sets $A_{n,2}$ of type 
{\bf II}$^{\prime}$
consisting of exactly 2 intervals are extremal ${\cal L}$-avoiding
sets in $[1,n]$ up to an error of at most $O(\log n)$ for every $n$. In
particular, these sets give the right value of $\lambda_{0}$. They 
do not attempt any classification of the extremal sets, however. They
also note that, whenever $c > (b+1)^{3/2}$, the same sets are of maximum size 
among all type {\bf II}$^{\prime}$ sets, and 
conjecture that they are still extremal, up to the same $O(\log n)$ error. 
\\
For {\sc Family II} equations they simply note that when $c > (2b)^{3/2}$, then
among all type {\bf II}$^{\prime}$ 
sets the largest consist of three intervals. They do not discuss whether such 
sets are extremal or not.
\\
\\
Our results concern the same two families of equations. For {\sc Family I}
we employ the methods of [1] to obtain a classification (Theorem 2.5) of
the extremal $\cal L$-avoiding sets whenever 
\begin{equation}
c > {(b+1)b^{2} \over b-1}.
\end{equation}
We show that, for every $n >>_{b,c} 0$ the sets $A_{n,2}$ are 
actually extremal and there are only a bounded number
of possibilities for the extremal subsets of $[1,n]$, all of which have 
a symmetric difference with $A_{n,2}$ of bounded size. Both bounds are 
independent of $n,b$ and $c$.
\par We show by means of an example that the lower bound (3) on $c$ cannot
be significantly improved, which also disproves the conjecture of 
Dilcher and Lucht. Namely we show that when $c = b^{2}$ another type of 
$\cal L$-avoiding subset of $[1,n]$ is larger than $A_{n,2}$ by a factor 
of $\Omega(n)$. In some cases we can prove that these other sets
are in fact extremal and conjecture that this is generally the case 
(conjecture 2.7).  
\\
\\
For {\sc Family II} equations we describe 
extremal sets in $[1,n]$ for all $n$,
and for every $b,c$ with $b > 1$ (Theorem 3.1). 
Their appearance takes three different forms, for values of $c$ in the 
following three ranges : (i) $c > 2b$, (ii) $2 \leq c < 2b$, (iii) $c = 1$. 
In contrast to when $b=1$, it is not the case for $c >> b$ that the
extremal sets consist of three intervals. Rather they 
are a hybrid between the 
two alternatives offered by our earlier Question.        

\section{Results for Family I equations} 

The methods of this section follow very closely those of [1], so we will 
not include full proofs of all results. Nevertheless, several technical
difficulties arise and considerable care is needed to dispose of them. 
Thus we will give a fair amount of detail anyway, even if the 
resulting computations become somewhat long-winded. 
Let $\cal L$ be a fixed equation
of the form $x+by = cz$ such that $b > 1$ and 
(3) holds. We prove analogues of Lemmas 
2,3,4 and Theorems 1,2 in [1]. First a definition corresponding to 
Definition 1 of that paper :
\\
\\
{\sc Definition 1} : Let $n \in {\hbox{{\bf N}}}$ and $A \subseteq [1,n]$
be $\cal L$-avoiding with smallest element $s := s_{A}$. Define sequences 
$(r_{i}), (l_{i}), (A_{i})$ by 
\begin{eqnarray*}
A_{0} := A, \;\;\; r_{1} := n, \\
l_{i} := \lfloor {b+1 \over c} r_{i} \rfloor, \;\;\; r_{i+1} := 
\lfloor {l_{i} + bs \over c} \rfloor, \\
A_{i} := ( A_{i-1} \backslash (r_{i+1},l_{i}] ) \cup [l_{i},r_{i}] \cap (s,n], 
\;\; {\hbox{for $i \geq 1$}}.
\end{eqnarray*}
Let $t$ denote the least integer such that $r_{t+1} < s$. Observe that for all 
$i \geq t$, 
\begin{equation}
A_{i} = A_{t} = [\alpha,r_{t}] \cup \left( \cup_{j=1}^{t-1} (l_{j},r_{j}] \right),
\end{equation}
where $\alpha = \alpha_{A} := \max \{l_{t}+1,s\}$. 
\\
\\
It is easy to see that, by construction, each set in the sequence $(A_{i})$ is
$\cal L$-avoiding provided $A_{0}$ is, and that $A_{t}$ is then an 
$\cal L$-avoiding set of type {\bf II}$^{\prime}$ in the introduction. 
The crucial
observation is the direct generalisation of Lemma 2(b) in [1] : because of its
importance and because an apparently awkward technicality arises in 
dealing with one of the cases ({\em Case I} below), 
we will provide a complete proof.
\\
\\
{\bf Lemma 2.1} {\em Let $n > 0$ and $A := A_{0} 
\subseteq [1,n]$ be $\cal L$-avoiding. 
Then $|A_{i}| \geq |A_{i-1}|$ for every $i > 0$.}
\\
\\
{\sc Proof} : Following the same reasoning as in [1], it suffices to 
prove the claim for $i = 1$, and thus to prove that, for every $n > 0$ and
every $\cal L$-avoiding subset of $[1,n]$, we have 
\begin{eqnarray}
|A| \leq |A \cap [1,r_{2,A}]| + \lceil \left(1-{b+1 \over c}\right) n 
\rceil,
\end{eqnarray}
where
\begin{eqnarray*}
r_{2,A} := \lfloor {\lfloor {b+1 \over c} n \rfloor + bs_{A} \over c} \rfloor.
\end{eqnarray*}
The proof is by induction on $n$, the case $n = 1$ being trivial. So suppose
the result holds for $1 \leq m < n$ and let $A$ be an $\cal L$-avoiding subset
of $[1,n]$. The result is again trivial if $s_{A} > \left( {b+1 \over c} 
\right) n$, so we may assume that $s_{A} \leq \left( {b+1 \over c} \right) n$
and thus that 
\begin{eqnarray*}
r_{2,A} \leq {\left( {b+1 \over c} \right) n + bs_{A} \over a} \leq 
\left( {b+1 \over c} \right)^{2} n < {n \over c},
\end{eqnarray*}
because of (3). 
\\
\\
First suppose that there exists $z \in A \cap \left( {n \over c}, 
{(b+1)n \over c} \right]$. To simplify notation,
denote $A^{c} := [1,n] \backslash A$.  
\\
\\
{\em Case I} : $z \leq bn/c$.
\\
\\
In this case we will show independently of the induction hypothesis that
something stronger than (5) holds, namely that 
\begin{equation}
|A| \leq \lceil \left( 1 - {b+1 \over c} \right) n \rceil.
\end{equation}
We have $cz \in (n,bn]$. Set $cz := t$. Now $A$ contains no solutions to the 
equation 
\begin{eqnarray}
x+by = t.
\end{eqnarray}
Hence for every $y \in [\frac{t-n}{b},\frac{t-1}{b}]$ at most one of the 
numbers $y$
and $t-by$ lies in $A$. Now $t-by \equiv t \; ({\hbox{mod $b$}})$ 
for every $y$. In this way we can locate in
$A^{c}$ at least as many numbers as there are 
numbers in the interval $[\frac{t-n}{b},\frac{t-1}{b}]$ not congruent to
$t \; ({\hbox{mod $b$}})$. Define two parameters $u,v \in [1,b]$ as 
follows :
\par (i) $u \equiv t \; ({\hbox{mod b}})$, 
\par (ii) the total number of integers in the interval
$\left[ \frac{t-n}{b}, \frac{t-1}{b} \right]$ is congruent to
$v \; ({\hbox{mod $b$}})$. 
\\
\\
Then one readily checks that the number of integers in $\left[ \frac{t-n}{b},
\frac{t-1}{b} \right]$ not congruent to $t \; ({\hbox{mod $b$}})$ is at
least
\begin{eqnarray*}
{b-1 \over b} \left[ {n-(u-1) \over b} - {b-v \over b-1} \right] = 
n \left( {b-1 \over b^{2}} \right) - {(b-1)(u-1) + b(b-v) \over b^{2}},
\end{eqnarray*}
and is at least one more than this when $v < b$ unless one of the first $v$ 
numbers in
the interval $\left[ \frac{t-n}{b}, \frac{t-1}{b} \right]$ is congruent to 
$t \; ({\hbox{mod $b$}})$. Set 
\begin{eqnarray*}
f(n,b,c,u,v) := {(b-1)(u-1) + b(b-v) \over b^{2}} - n \left(
{b-1 \over b^{2}} - {b+1 \over c} \right).
\end{eqnarray*}
Note that (3) implies that $f(n,b,c,u,v) < 2$ but, for (6) to be
already satisfied we would need $f(n,b,c,u,v) < 1$. This is where things 
get messy. Note that certainly $f(n,b,c,u,v) < 1$ unless perhaps
$u \geq 3$, $v \leq u-2$ and one of the first $v$ numbers in the interval
$\left[ \frac{t-n}{b}, \frac{t-1}{b} \right]$ is 
congruent to $u \; ({\hbox{mod $b$}})$. The first assumption implies
in particular that $b \geq 3$. All three together imply that the 
numbers 1 and 2 both lie to the left of the interval $\left[ 
\frac{t-n}{b}, \frac{t-1}{b} \right]$, and neither is congruent to
$t \; ({\hbox{mod $b$}})$. Thus neither can have been already located in 
$A^{c}$ via the pairing arising from (7). Since it suffices at this point to
locate just one extra element in $A^{c}$ we may for the remainder of
this argument assume that $b \geq 3$ and that $1,2 \in A$. The latter implies
that there are no solutions in $A$ to either of the equations 
\begin{eqnarray}
x+by = c, \\ x+by = 2c.
\end{eqnarray}
To continue the argument we go back to (7). To locate elements in 
$A^{c}$ we paired off numbers in $\left[ \frac{t-n}{b}, \frac{t-1}{b} \right]$
not congruent to $t \ ({\hbox{mod $b$}})$ with numbers in $[1,n]$
congruent to $t \; ({\hbox{mod $b$}})$. It would thus suffice if
we could also pair off at least one further number in the former interval,
whch we call $\cal I$.
It is easy to see that this can definitely be done if $\cal I$ contains
a total of at least $b^{2}$ numbers. Hence we may further assume now that
\begin{equation}
| {\cal I} | \leq b^{2},
\end{equation}
and hence that $n \leq b^{3}$, though we will make no explicit use of this 
latter fact. 
\par From (10) we 
want to conclude that either $n < c$ or, 
for an appropriate choice of the original $z$, that $c \equiv t \equiv 
0 \; ({\hbox{mod $b$}})$. So suppose $n \geq c$. First set
\begin{eqnarray*}  
x_{1} := 1, \;\; x_{2} := b+1, \;\; y_{1} := c-b, \;\; y_{2} := c - b(b+1).
\end{eqnarray*}
Since $A$ contains no solutions to (8), at most one of 
$x_{i}$ and $y_{i}$ is in $A$ for each $i=1,2$. Thus $y_{1} \in A^{c}$ since 
$n \geq c$ and we 
already know that $x_{1} \in A$. From (3) it follows that
$x_{2} < y_{2}$ and also from (10) that $y_{1} - y_{2} > |{\cal I}|$, so
that at least one of $y_{1}$ and $y_{2}$ must lie outside $\cal I$. 
Furthermore, since $u \geq 3$, neither of the $x_{i}$ is congruent to 
$t \; ({\hbox{mod $b$}})$. If $b+1 \not\in {\cal I}$ it is thus
already clear that, unless $c \equiv t \; ({\hbox{mod $b$}})$, we can find
amongst $x_{2},y_{1},y_{2}$ at least one element of $A^{c}$ not previously
located via (7). But similarly, if $b+1 \in {\cal I}$ then one easily
checks that (3) implies that $y_{1} \not\in {\cal I}$ and so we have the same
conclusion. 
\par Thus we are done if $n \geq c$ unless $c \equiv t \; ({\hbox{mod $b$}})$.
To get that $c \equiv 0 \; ({\hbox{mod $b$}})$ it would then suffice to
also show that $2c \equiv t \; ({\hbox{mod $b$}})$. If $n \geq 2c-b$ then
this is immediately 
achieved by a similar argument to the one just given, but this time using (9)
instead of (8). If $n < 2c-b$ then we just have to note that we could
have from the beginning chosen $z := 2$, in which case $2c = t$, by definition.
\\
\\
Thus (6) holds unless either $n < c$ or $c \equiv 0 \; ({\hbox{mod $b$}})$.
By (3) the latter would imply that $c \geq b^{2} + kb$ where we can take
$k = 3$ when $b > 3$ and $k=4$ when $b=3$. Then
\begin{eqnarray*}
f(n,b,c,u,v) = {(b-1)(u-1) + b(b-v) \over b^{2}} - n \left( {b-1 \over b^{2}}
- {b+1 \over c} \right) \\
\leq {(b-1)(2b-1) \over b^{2}} - n \left( {b-1 \over b^{2}} - 
{b+1 \over b^{2}+kb} \right).
\end{eqnarray*}
We'll be done unless $f(n,b,c,u,v) \geq 1$. One checks that this already
forces $n < c$ when $b = 3$, and that for $b > 3$ it
yields (taking $k=3$) that
\begin{eqnarray*}
n \leq b^{2}+3b+1+{6 \over b-3} \leq c+1+{6 \over b-3}.
\end{eqnarray*}
This in turn yields that 
\begin{equation}
\lfloor \left( {b+1 \over c} \right) n \rfloor
\leq b+1
\end{equation} 
except when $b=4, c=28$ and $n \in \{34,35\}$. One 
easily checks that (6) holds in these exceptional cases, so we're left with 
(11). First suppose $n \geq c$ so that 
$\lfloor \left({b+1 \over c} \right) n \rfloor = b+1$, and consider (8). 
For $1 \leq i \leq b+1$ set $y_{i} := i$ and $x_{i} := c - ib$. Then
at least one of $x_{i}$ and $y_{i}$ is in $A^{c}$ for each $i$. But (3)
implies that $x_{b+1} < y_{b+1}$, hence $|A^{c}| \geq b+1$ which proves (6).
\par We are thus indeed 
left with the case when $n < c$. Now we could have chosen 
$z := 1$ initially and thus paired off numbers in ${\cal I}_{1} := 
\left[ \frac{c-n}{b}, 
\frac{c-1}{b} \right]$ not congruent to $c \; ({\hbox{mod $b$}})$ with numbers
in $[1,n]$ congruent to $c \; ({\hbox{mod $b$}})$. As usual, it 
suffices to locate
at least one further element in $A^{c}$. First suppose 
\begin{equation}
{2c -1 \over b} \leq n.
\end{equation}
Then similarly, by (9), we can pair off numbers in ${\cal I}_{2} := 
\left[ \frac{2c-n}{b}, 
\frac{2c-1}{b} \right]$ not congruent to $2c \; ({\hbox{mod $b$}})$ with
numbers in $[1,n]$ congruent to $2c \; ({\hbox{mod $b$}})$. The crucial 
point is that, since $n < c$, the intervals ${\cal I}_{1}$ and ${\cal I}_{2}$
are disjoint. Each interval certainly contains at least three elements by (12).
It is then easy to see that the ${\cal I}_{2}$-pairing will certainly locate
at least one more element in $A^{c}$ unless, at the very least, 
$2c \equiv c \equiv 0 \; ({\hbox{mod $b$}})$. But in that case the map
$\phi : y \mapsto y + \frac{c}{b}$ is a bijection from ${\cal I}_{1}$
to ${\cal I}_{2}$ so that if the ${\cal I}_{1}$-pairing pairs $y$ with
$x$, say, then the ${\cal I}_{2}$-pairing pairs $\phi(y)$ with $x$. 
If we now choose $y$ as the smallest multiple of $b$ in ${\cal I}_{1}$, then
we see that one of the two pairings must locate the desired extra element in
$A^{c}$, unless perhaps $\frac{c}{b} \equiv 0 \; ({\hbox{mod $b$}})$ also. 
But then $c \equiv 0 \; ({\hbox{mod $b^{2}$}})$ and thus $c \geq 2b^{2}$ if
$b > 3$ and $c \geq 3b^{2}$ if $b = 3$. But then, calculating as before, 
we'll have $f(n,b,c,u,v) < 1$ unless perhaps 
\begin{eqnarray*}
\left\{ \begin{array}{lr} n \leq {2(b^{2}-3b+1) \over b-3}, & {\hbox{when
$b > 3$}}, \\ n \leq {3(b^{2}-3b+1) \over 2(b-2)}, & {\hbox{when $b=3$}}.
\end{array} \right.
\end{eqnarray*}
But these inequalities contradict (12). Now we are only left with the 
possibility that $n < \frac{2c-1}{b}$, hence that $\lfloor 
\left( {b+1 \over c} \right) n \rfloor \in \{1,2\}$. But in each case one 
may check that one can 
locate one or two elements of $A^{c}$ as approppriate, by considering 
solutions of (8) with $x$ and $y$ close to $\frac{c}{b+1}$. This finally
completes the analysis of {\em Case I}.  
\\
\\
{\em Case II} : $z > bn/c$.
\\
\par Then $cz \in (bn,(b+1)n]$. Put $cz := t$ again. Let $t = (b+1)n-s$ where 
$0 \leq s < n$. If $x+by = t$ for some integers $x,y \in [1,n]$, then 
$y \geq n - s/b$ and $x \geq n-s$. Since $A$ avoids $\cal L$ 
we thus find, for every 
integer $y \in A \cap [n - \frac{s}{b},n]$ not congruent to $t \; 
({\hbox{mod $b$}})$, 
an integer
$x \in A^{c} \cap [n-s,n]$, congruent to $t \; ({\hbox{mod $b$}})$. Noting in 
addition that at least one of $n$ and $n-s$ is not in $A$, one 
readily verifies that hence 
\begin{equation}
|A \cap [n-s,n]| \leq \left( 1 - {b-1 \over b^{2}} \right) (s+1).
\end{equation}
We now apply the induction hypothesis. Let $B := A \cap [1,n-s-1]$. 
If $B$ is empty then (13) and (3) immediately imply (5). Otherwise
clearly $s_{B} = s_{A}$ and $r_{2,B} \leq r_{2,A}$, so the induction
hyothesis gives that
\begin{eqnarray*}
|A| \leq |A \cap [1,r_{2,A}]|
+ \lceil \left( 1 - {b+1 \over c} \right) (n-s-1) \rceil + \left( 1 - 
{b-1 \over b^{2}} \right) (s+1), 
\end{eqnarray*}
from which (5) follows by another application of (3).
\\
\\
We have thus completed the induction step under the assumption that 
$A \cap \left( {n \over c}, {(b+1)n \over c} \right] \neq \phi$, so we
can now assume the intersection is empty. Suppose $z \in A \cap (r_{2,A},n/c]$.
Then $\lfloor {b+1 \over c} n \rfloor + bs_{A} < cz \leq n$ and 
$cz - bs_{A} \in A^{c}$. In other words, we can pair off
elements in $A \cap (r_{2,A},\frac{n}{c}]$ with elements in 
$\left( {b+1 \over c} n, n\right] \cap A^{c}$. This immediately implies (5)
and completes the proof of Lemma 2.1.  
\\
\\
{\bf Lemma 2.2} {\em Let $A$ be an $\cal L$-avoiding subset of $[1,n]$ of
maximum size.
Let $s = s_{A}$ and $t := \max \{i \in {\hbox{{\bf N}}} : r_{i} \geq s\}$. If
$n >>_{b,c} 0$ then $t = 2$.}
\\
\\
{\sc Proof} : Just follow the reasoning in the proof of Lemma 3 in [1]. By
Lemma 2.1, it suffices to 
know that there exists an absolute positive constant $\kappa^{0}_{b,c}$
such that, if $t \neq 2$ then 
\begin{eqnarray*}
{|A| \over n} \leq D(b,c) - \kappa^{0}_{b,c},
\end{eqnarray*}
where 
\begin{eqnarray}
D(b,c) := {(c-b-1)(c^{2}-b^{2}+1) \over c[c^{2}-b(b+1)] }
\end{eqnarray}
is such that, in the notation of eq.(4), $|A_{2}| = D(b,c) \cdot n + O(1)$ when
$s = l_{2}+1$. The core of a proof that such a constant exists is contained 
in the proof of Lemma 1 in [4], though one has to be a little careful
since there only sets $A_{t}$ in which $s = l_{t}+1$ are considered. 
However one can tediously check that allowing for arbitrary $s \in 
(l_{t},r_{t}]$ will not change matters (I note that the authors of [4]
needed this fact in Section 3 of their paper, though they do not seem to
explicitly mention it anywhere).
\\
\\
{\bf Lemma 2.3} {\em Let $n >>_{b,c} 0$. If $A$ is an $\cal L$-avoiding
subset of $[1,n]$ of maximum size then there exists an absolute positive
constant $\kappa^{1}_{b,c}$ such that $S-\kappa^{1}_{b,c} \leq s_{A} \leq
S+2$ where $S = \lfloor {(b+1)^{2} n \over c[c^{2}-b(b+1)]} \rfloor$.}
\\
\\
{\sc Proof} : The proof follows that of Lemma 4 in [1]. We set 
\begin{eqnarray*}
s^{\prime} := \min \{s \in [1,n] : l_{2}(s) < s \}.
\end{eqnarray*}
A computation similar to that in the Appendix of [1] yields that
\begin{eqnarray*}
l_{2}(s) < s \Leftrightarrow s > {(b+1)^{2} \over c[c^{2}-b(b+1)]} n - 
(\epsilon_{1}c + \epsilon_{2})\left({b+1 \over c^{2}-b(b+1)}\right),
\end{eqnarray*}
where $\epsilon_{1}, \epsilon_{2} \in [0,1)$. By (3) it follows that
\begin{equation}
s^{\prime} \in [S,S+1].
\end{equation}
Now we have 
\begin{eqnarray*}
|A_{2}(s)| = \left\{ \begin{array}{lr} \lceil \left( 1 - {b+1 \over c}
\right) n \rceil + r_{2}(s) - s +1, & {\hbox{if $s \geq s^{\prime}$}}, \\
\lceil \left( 1 - {b+1 \over c} \right) n \rceil + r_{2}(s) - l_{2}(s), & 
{\hbox{if $s < s^{\prime}$}}. \end{array} \right.
\end{eqnarray*}
First suppose $s \geq s^{\prime}$. We will certainly have 
$|A_{2}(s^{\prime}+2)| < |A_{2}(s)|$ because of (3) and 
since $r_{2}(s)$ can only increase at
most once in $\lfloor \frac{c}{b} \rfloor$ times. This proves that 
$s_{A} \leq S+2$ 
for a maximum $\cal L$-avoiding $A$. Secondly, if $s < s^{\prime}$
then $|A_{2}(s)|$ will decrease once $r_{2}(s) - l_{2}(s)$ decreases, and
clearly this will happen after $\Omega(\frac{c}{b})$ steps. This proves
that $S - \kappa^{1}_{b,c} \leq s_{A}$ for a maximum $A$ and some 
$\kappa^{1}_{b,c} = \Omega(\frac{c}{b})$. 
\\
\\
{\bf Theorem 2.4} {\em $r_{\cal L}(n) = D(b,c) \cdot n + O(1)$, where $D(b,c)$
is given by (14). In particular, $\lambda_{0,{\cal L}} = D(b,c)$.}
\\
\\
{\sc Proof} : See Lemma 1 in [4] and Theorem 1 in [1]. 
Note that the second statement is 
the same as in Theorem 1 of [4], 
just with a better lower bound for $c$, namely (3).
\\
\\
We can now present the main classification result, analogous to Theorem 2 in 
[1]. In fact the result here is in some sense even cleaner, as the 
maximum $\cal L$-avoiding sets consist essentially of two rather than three
intervals and there is thus even less possibility for variation. 
\\
\\
{\bf Theorem 2.5} {\em Let $b,c \in {\hbox{{\bf N}}}$ with $b > 1$ and 
$c$ satisfying (3). Let $\cal L$ be the equation $x+by=cz$. 
Let $n > 0$. Define $S$ and $s^{\prime}$ as in Lemma 
2.3. Let $A$ be an $\cal L$-avoiding subset of $[1,n]$ of maximum 
size and with smallest element $s = s_{A}$. 
If $n >>_{b,c} 0$ then the following holds : $s \in [S,S+2]$
and $A = {\cal I}_{2} \cup {\cal I}_{1}$ where 
\begin{eqnarray}
{\cal I}_{2} \in \left\{ \begin{array}{lr} \{[s,r_{2}], [s,r_{2}+1]\}, & 
{\hbox{if $s \geq s^{\prime}$}}, \\ \{[s,r_{2}), [s,r_{2}]\backslash\{r_{2}-
\xi_{1}\}\},
& {\hbox{if $s < s^{\prime}$}}, \end{array} \right.
\end{eqnarray}
for some $\xi_{1} \in [1,b]$, and 
\begin{eqnarray}
{\cal I}_{1} \in \left\{ \begin{array}{lr} \left.(l_{1},n]\right., & 
{\hbox{if $r_{2}+1 \not\in A$ and $l_{1} \not\in A$}}, \\ 
\left.[l_{1},n] \backslash \{n-\xi_{2}\}\right.,
& {\hbox{if $r_{2}+1 \not\in A$ and $l_{1} \in A$}}, \\
\left.(l_{1},n]\right. \backslash \{l_{1} + \xi_{3}\}, 
& {\hbox{if $r_{2}+1 \in A$ and $l_{1} \not\in A$}}, \\ 
\left.[l_{1},n] \backslash \{l_{1}+\xi_{4},n-\xi_{5}\}\right., 
& {\hbox{if $r_{2}+1 \in A$ and $l_{1} \in A$}}, \end{array} \right.
\end{eqnarray}
for some $\xi_{2} \in [0,n]$, $\xi_{3} \in [1,n]$, $\xi_{4} \in [1,b-1]$, 
$\xi_{5} \in [0,n]$. Moreover, in all cases where they arise, the parameters
$\xi_{i}$, $i = 1,...,5$, are uniquely determined by $n$ and $s$ 
according to the following relations :
\begin{eqnarray}
cl_{2} = (b+1)r_{2} - \xi_{1}, \\
cl_{1} = (b+1)n - \xi_{i}, \;\; i \in \{2,5\}, \\
c(r_{2}+1) = bs + (l_{1} + \xi_{i}), \;\; i \in \{3,4\}.
\end{eqnarray}}
\\
\\
{\sc Remark} : As is the case in [1], Theorem 2.5 does not precsiely determine
the maximum $\cal L$-avoiding subsets of $[1,n]$ for every $n >>_{b,c} 0$.
For any particular $n$, some of the possibilities listed for $A$ may either 
not be $\cal L$-avoiding or not have maximum size. But the important point
is that we have a bounded number of possibilities for any $n$, and the 
symmetric difference between any two of these possibilities is also 
bounded in size, both bounds being independent of $n,b$ and $c$. Since the
periodicity phenomenon described in Section 4 of [1] easily
generalises to the present setting, Theorem 2.5 thus reduces the 
precise classification of the extremal $\cal L$-avoiding sets to a 
finite computation for any given $\cal L$. 
\\
\\ 
{\sc Proof of Theorem 2.5} : Again we follow the approach in [1]. On the one 
hand, since $t=2$ here, there are fewer steps in the analysis. On the 
other hand, we will need a somewhat modified argument in one of the steps.
We shall thus present the full argument quite carefully, though not in 
every single detail. We will need something analogous to Lemma 1 of [1]. 
Let $z \in {\hbox{{\bf N}}}$. Then $x+by=cz$ will have solutions where $x$ and
$y$ are close to $\frac{c}{b+1}$. Indeed if $1-\epsilon := \frac{c}{b+1} - 
\lfloor \frac{c}{b+1} \rfloor$ then we have a solution $(x_{0},y_{0}) = 
\left( \frac{c}{b+1} -b\epsilon, \frac{c}{b+1} + \epsilon \right)$. For each 
$i \geq 0$ define the solution
\begin{eqnarray*}
(x_{i},y_{i}) := \left( \frac{c}{b+1} - b(i + \epsilon), \frac{c}{b+1} + 
i + \epsilon \right).
\end{eqnarray*}
Now for any $z \in {\hbox{{\bf N}}}$ and $d \geq 0$ let $I^{d}_{z}$ denote
the interval $[x_{d},y_{d}]$. Then the analogue of Lemma 1 we need is 
\\
\\
{\bf Lemma 2.6} {\em Let $A \subseteq [1,n]$ be $\cal L$-avoiding and 
$z \in A$. Then for any $d \geq 0$, $|I^{d}_{z} \backslash A| \geq d+1$.}
\\
\\
Let $A$ be a maximum $\cal L$-avoiding subset of $[1,n]$. By Lemma 2.1
we know that $|A| = |A_{1}| = |A_{2}|$ and by Lemma 2.2 that $r_{3} < s$, when
$n$ is sufficiently large. 
The proof of Theorem 2.5 is accomplished in three steps. First, by 
comparing $A_{2}$ with $A_{1}$ we show that $A$ contains almost the whole 
interval $(l_{2},r_{2}]$. Here the argument entirely parallells that in [1].
In the second step we deduce that $(r_{2},l_{1}] \cap A$ is almost empty.
Here we use Lemma 2.6, but in a somewhat different way than in [1], as we
will instead use an approach similar to that in the proof of Lemma 2.1. 
The final step is to compare $A$ with $A_{1}$ to show that $A$ contains almost
all of $(l_{1},n]$. 
\\
\\
{\em Step 1} : If $ s > l_{2}$ then clearly $[s,r_{2}] \subseteq A$. Lemma 2.2
and (15) give $s \in [S,S+2]$. Suppose $s \leq l_{2}$. We want to show that
$s = l_{2}$. Suppose on the contrary that $B := [s,l_{2})$ is non-empty. 
Put 
\begin{eqnarray*}
C := I^{1}_{s} \cup \bigcup_{b \in B \backslash \{s\}} I^{0}_{b}.
\end{eqnarray*}
It is clear that $C \subseteq [l_{2},r_{2}]$ for all $n >>_{b,c} 0$. The 
crucial point is that (3) guarantees that all the intervals making up $C$
are pairwise disjoint. Thus Lemma 2.6 implies that $|C \backslash A| > |B|$,
which contradicts the maximality of $A$. Thus $s = l_{2}$, which implies
on the one hand (computation required, using (3)) that $s \geq S$, and on the
other that 
\begin{equation}
|A \cap [s,r_{2}]| = |(s,r_{2}]|.
\end{equation}
If $r_{2} \not\in A$ we infer that $A \cap [s,r_{2}] = [s,r_{2})$. If 
$r_{2} \in A$ then $cs -br_{2} = cl_{2} -br_{2} \not\in A$, so $-c+1 \leq 
cl_{2}-(b+1)r_{2} \leq -1$ and hence 
\begin{eqnarray*}
{-c+1 \over b+1} \leq {c \over b+1} l_{2} - r_{2} < 0.
\end{eqnarray*}
But if ${c \over b+1} l_{2} - r_{2} \leq -1$ then $I^{1}_{l_{2}} \subseteq
(l_{2},r_{2}]$ (for $n >>_{b,c} 0$), which would 
contradict (21). Thus ${c \over b+1} l_{2} - r_{2} \in (-1,0)$, which confirms
that $A \cap [s,r_{2}] = [s,r_{2}] \backslash \{r_{2} - \xi_{1}\}$ for the 
unique $\xi_{1} \in [1,b]$ satisfying (18). We also 
deduce that $r_{2}+1 \not\in A$ since $cl_{2} - b(r_{2}+1) < r_{2} - b$ and thus
lies in $A$. 
\par This completes Step 1 and shows that $A$ must contain a set ${\cal I}_{2}$
which is one of the possibilities given by (16). 
\\
\\
{\em Step 2} : We will show that if $n >>_{b,c} 0$ then $A \cap
[r_{2}+2,l_{1}) = \phi$. We have 
\begin{eqnarray*}
|A| = |A_{1}| = |(A \backslash (r_{2},l_{1}]) \cup (l_{1},n]|.
\end{eqnarray*}
So the idea is to show that if $A \cap [r_{2}+2,l_{1})$ were non-empty, then 
it would have to have smaller cardinality than $(l_{1},n] \backslash A$. 
Since $b > 1$ the direct analogue of the
argument in [1] will not work. Instead we gain inspiration from
the proof of Lemma 2.1. First suppose there exists $z \in A \cap 
\left( {n \over c},
{bn \over c} \right]$. As in the proof of Lemma 2.1, this implies that
\begin{eqnarray*}
|A| < \left( 1 - {b-1 \over b^{2}} \right) n + 2,
\end{eqnarray*}
and thus, by Theorem 2.4, $A$ can't possibly be maximum $\cal L$-avoiding for
$n >>_{b,c} 0$. Next suppose there exists $z \in A \cap \left( {bn \over c},
{(b+1)n \over c} \right]$. Again, as in the proof of Lemma 2.1, this will
gives us a $\sigma \in [0,n]$ such that 
\begin{eqnarray*}
|A \cap [n-\sigma,n]| \leq \left( 1 - {b-1 \over b^{2}} \right) (\sigma +1).
\end{eqnarray*}
Clearly then, by Theorem 2.4, there exists a constant $\kappa^{2}_{b,c} > 0$
such that $A$ cannot possibly be maximum $\cal L$-avoiding if 
$\sigma > \kappa^{2}_{b,c}$. 
Thus there exists a corresponding $\kappa^{3}_{b,c} > 0$
such that we can now deduce that $A \cap \left( {n \over c}, l_{1} - 
\kappa^{3}_{b,c} \right]$ is empty. 
\par Let ${\cal U}_{1} := A \cap [r_{2}+2,n/c]$ and 
${\cal U}_{2} := A \cap (l_{1} - \kappa^{3}_{b,c},l_{1})$. It 
remains to show that ${\cal U}_{1}$ and ${\cal U}_{2}$ are empty, so let us
suppose otherwise. 
\par For 
$z \in {\cal U}_{1}$, let ${\cal C}_{z} := \{cz-bs,cz-b(s+1)\}$.
Also let $C_{r_{2}+1} := \{c(r_{2}+1) - bs\}$. 
Then ${\cal C}_{z} \cap A = \phi$ for any $z \in A \cap 
[r_{2}+1,n/c]$. Clearly, if $n >>_{b,c} 0$ then 
$C_{z} \subset (l_{1},n - \Omega(n)]$. Also (3) guarantees that the 
${\cal C}_{z}$ are pairwise disjoint. 
\par For $z \in {\cal U}_{2}$ let ${\cal D}_{z} := I^{1}_{z} \cap A^{c}$ if
$z$ is the smallest element of ${\cal U}_{2}$ and ${\cal D}_{z} := 
I^{0}_{z} \cap A^{c}$ otherwise. Again (3) guarantees that the ${\cal D}_{z}$
are pairwise disjoint. Clearly there exists a constant $\kappa^{4}_{b,c} > 0$
such that all the ${\cal D}_{z}$ are contained in $[n - \kappa^{4}_{b,c},n]$.
Thus the ${\cal D}_{z}$ are also disjoint from the ${\cal C}_{z}$. 
\par In summary, we can thus conclude that, for $n$ sufficiently large,
\begin{eqnarray*}
|(l_{1},n] \cap A^{c}| \geq \delta_{1} + 2|{\cal U}_{1}| + (|{\cal U}_{2} + 1 - 
\delta_{2}),
\end{eqnarray*}
where $\delta_{1} = 1$ if $r_{2} + 1 \in A$ and zero otherwise, and
$\delta_{2} = 1$ if ${\cal U}_{2} = \phi$ and zero otherwise. It follows
immediately that $|A| < |A_{1}|$ unless ${\cal U}_{1}$ and ${\cal U}_{2}$
are both empty. This completes Step 2.
\\
\\
{\em Step 3} : It just remains to show that the possibilities for $A \cap
[l_{1},n]$ are as given by (17). By Steps 1 and 2 we only have four cases to 
consider, according to $A \cap \{r_{2}+1,l_{1}\}$. To verify the 
various possibilities for ${\cal I}_{1}$ one considers the numbers
$c(r_{2}+1) - bs$ or $cl_{1} - (b+1)n$ as appropriate, the analysis being 
similar to that in the latter part of Step 1 above. We omit the details
and consider the proof of Theorem 2.5 as complete.        
\\
\\
We close this section by showing that, in general, the bound (3) cannot be 
significantly decreased without the extremal sets avoiding 
$x+by=cz$ looking quite different than those described in the above
theorem. In
particular we have a counterexample to the conjecture in [4] that a 
bound of $c > (b+1)^{3/2}$ should suffice. For a counterexample we
set $c = b^{2}$. Let ${\cal L}_{b}$ denote the equation 
$x+by=b^{2}z$ where $b > 1$. The constructions considered earlier in this
section yield that
\begin{eqnarray*}
\lambda_{0,{\cal L}_{b}} \geq D(b,b) = {(b^{2}-b-1)(b^{4}-b^{2}+1) \over 
b^{2}[b^{4} - b(b+1)]}.
\end{eqnarray*}
But for every $b > 1$ the true value of $\lambda_{0,{\cal L}_{b}}$ is larger.
For let 
\begin{eqnarray*}
A_{b} := \{ u \cdot b^{3i} : u > 0, \; i \geq 0 \; {\hbox{and $b \dagger u$}}
\}.
\end{eqnarray*}
Then clearly $A_{b}$ is an ${\cal L}_{b}$-avoiding subset of {\bf N}
and 
\begin{eqnarray*}
d(A_{b}) = {b^{2} \over b^{2}+b+1} > {(b^{2}-b-1)(b^{4}-b^{2}+1) \over
b^{2}[b^{4}-b(b+1)]} \;\;\; \forall \; b \geq 2.
\end{eqnarray*}
We conjecture the following :
\\
\\
{\bf Conjecture 2.7} {\em For every $n > 0$ and every $b \geq 2$ the 
set $A_{b} \cap [1,n]$ is an ${\cal L}_{b}$-avoiding subset of $[1,n]$ of
maximum size. In particular $\lambda_{0,{\cal L}_{b}} = \rho_{{\cal L}_{b}}
= {b^{2} \over b^{2}+b+1}$.}
\\
\\
We suspect in fact that for $n >>_{b} 0$ any extremal ${\cal L}_{b}$-avoiding
subset of $[1,n]$ must be very similar to $A_{b} \cap [1,n]$. 
Frustratingly we have not been able to verify any of these assertions
in general, not even the value of $\lambda_{0,{\cal L}_{b}}$. We do have
proofs of Conjecture 2.7 for $b = 2,3$ which we now present. They employ
the same idea, but things get pretty messy for $b=3$ and we don't 
see how to make the same idea work for larger $b$. 
\\
\\
{\bf Theorem 2.8} {\em Conjecture 2.7 holds for $b = 2$ and $b = 3$.}
\\
\\
{\sc Proof for $b=2$} : Fix $n > 0$. Put $A = A_{2} \cap [1,n]$. Let $B$ be an
${\cal L}_{2}$-avoiding subset of 
$[1,n]$. We must show that $\mid B \mid \; \leq \; \mid A \mid$. We will do 
this by exhibiting a one-to-one function 
\begin{eqnarray*}
f : B \backslash A \rightarrow A \backslash B.
\end{eqnarray*}
For $k=1,2$ let 
\begin{eqnarray*}
B_{k} := B \cap \{ u \cdot 2^{3i+k} : 2 \dagger u, i \geq 0\}.
\end{eqnarray*}
Then $B \backslash A = B_{1} \sqcup B_{2}$. For $x \in B \backslash A$, we shall define $f(x)$ according as to whether $x \in B_{1}$ or $x \in B_{2}$. 
\\
\\
First suppose $x \in B_{1}$. Then $x=2y$ for some $y \in A$. But $y \not\in B$ since $B$ avoids ${\cal L}_{2}$ and 
\begin{eqnarray*}
4 \cdot y = 2 \cdot y + (2y).
\end{eqnarray*}
So in this case we define $f(x) = y$. 
\\
\\
Next suppose $x \in B_{2}$. Then $x=4y$ for some $y \in A$. Then $3y \in A$, but $3y \not\in B$ since $B$ avoids ${\cal L}_{2}$ and 
\begin{eqnarray*}
4(3y) = 2(4y) + (4y).
\end{eqnarray*}
So in this case we define $f(x) = 3y$.
\\
\\
It is clear that the restrictions of $f$ to both $B_{1}$ and $B_{2}$ are one-to-one. So it remains to show that $f(B_{1}) \cap f(B_{2}) = \phi$. 
\par So let $y,z \in A$ and suppose that $f(2y) = f(4z)$. Thus $y=3z$ and so $2y=6z$. So both $4z \in B$ and $6z \in B$. But this is a contradiction, since $B$ avoids ${\cal L}_{2}$ and
\begin{eqnarray*}
4(4z) = 2(6z) + (4z).
\end{eqnarray*}
{\sc Proof for $b=3$} : Fix $n > 0$. Put $A = A_{3} \cap [1,n]$ and let $B$ be 
an ${\cal L}_{3}$-avoiding subset of 
$[1,n]$. As before we will describe an explicit one-to-one function 
$f : B \backslash A \rightarrow A \backslash B$. We define the sets $B_{1}$ 
and $B_{2}$ in an analogous manner to above and for $x \in B \backslash A$ 
will define $f(x)$ according as to whether $x \in B_{1}$ or $B_{2}$. This 
time, both the definition of $f$ and the proof that it is one-to-one will be 
somewhat more complicated than before, so we divide this process into three 
clear steps. 
\\
\\
{\em Step 1} : We define $f$ on $B_{1}$ and show that $f \mid_{B_{1}}$ is one-to-one.
\\
\\
Let $x \in B_{1}$. Then $x=3y$ for some $y \in A$. Then $2y \in A$. Now it can't be the case that both $y$ and $2y$ lie in $B$, since $B$ avoids
${\cal L}_{3}$ and 
\begin{eqnarray*}
9 \cdot y = 3(2y) + (3y).
\end{eqnarray*}
Hence we define 
\begin{eqnarray*}
f(x) = \left\{ \begin{array}{lr} y, & {\hbox{if $y \not\in B$}}, \\ 2y, & {\hbox{otherwise}}. \end{array} \right.
\end{eqnarray*}
We need to show that $f$ is on-to-one on $B_{1}$. Suppose otherwise. 
Then there is an $x=3y \in B_{1}$ such that $2x=6y \in B_{1}$ 
and $f(x)=f(2x)=2y$. But this implies that $y \in B$, which is a 
contradiction, since $B$ avoids ${\cal L}_{3}$ and 
\begin{eqnarray*}
9 \cdot y = 3 \cdot y + (6y).
\end{eqnarray*}
{\em Step 2} : We define $f$ on $B_{2}$ in such a way that 
\begin{eqnarray*}
f(B_{1}) \cap f(B_{2}) = \phi.
\end{eqnarray*}
Let $x \in B_{2}$. Then $x=9y$ for some $y \in A$. Now $4y \in A$, but 
$4y \not\in B$ since $B$ avoids ${\cal L}_{3}$ and 
\begin{eqnarray*}
9(4y) = 3(9y) + (9y).
\end{eqnarray*}
If $4y \in f(B_{1})$ then either 
\\
\par (a) $12y \in B$, in which case $f(12y)=4y$, or 
\par (b) $12y \not\in B$, but $6y \in B$ and $f(6y)=4y$. In this case, the definition of $f$ on $B_{1}$ implies that $2y \in B$. I also claim that in this case, $y \not\in B \cup f(B_{1})$. Suppose $y \in B$. Then, since $6y \in B$, 
the equation $9 \cdot y = 3 \cdot y + (6y)$ contradicts the fact that $B$ 
avoids ${\cal L}_{3}$. So suppose $y \in f(B_{1})$. 
Then either $y = f(3y)$ or $y = f(3y/2)$. But if $3y \in B$ then, 
since both $6y$ and $9y$ are also in $B$, the equation $9(3y)=3(6y)+(9y)$ 
contradicts the fact that $B$ avoids ${\cal L}_{3}$. 
And if $3y/2 \in B$ then the equation 
$9 (\frac{3y}{2} ) = 3(\frac{3y}{2}) + (9y)$ likewise gives a 
contradiction. 
\\
\\
Hence, we begin by defining 
\begin{eqnarray*}
f(x) = \left\{ \begin{array}{lr} 4y, & {\hbox{if $12y \not\in B$ and $\{2y,6y\} \not\subseteq B$}}, \\ y, & {\hbox{if $12y \not\in B$ and $\{2y,6y\} \subseteq B$}}. \end{array} \right.
\end{eqnarray*}
It remains to define $f(x)$ when $12y \in B$. Notice that then $2y \not\in B$, since $9(2y)=3(2y)+(12y)$. If $2y \in f(B_{1})$ then either $2y=f(3y)$ or $2y=f(6y)$. So we define 
\begin{eqnarray*}
f(x)=2y, \;\;\;\; {\hbox{if $12y \in B$, $3y \not\in B$ and $6y \not\in B$}}.
\end{eqnarray*}
Next suppose $3y \not\in B$ but $6y \in B$. Then $y \not\in B$ since $9\cdot y = 3 \cdot y + (6y)$. And $y \not\in f(B_{1})$ either, since if it were then either
$y=f(3y)$ or $y=f(3y/2)$. But $3y \not\in B$, by assumption, and $3y/2 \not\in B$ since $9 (\frac{3y}{2}) = 3(\frac{3y}{2}) + (9y)$. 
\par Thus we may define 
\begin{eqnarray*}
f(x) = y, \;\;\;\; {\hbox{if $12y \in B$, $3y \not\in B$ and $6y \in B$}}.
\end{eqnarray*}
Finally, it remains to define $f$ when both $12y$ and $3y$ are in $B$. 
I claim that in this case, $8y \not\in B \cup f(B_{1})$. Suppose $8y \in B$. 
Then $9(3y)=3(8y)+(3y)$, contradicting $B$'s avoidance of ${\cal L}_{3}$. 
Suppose $8y \in f(B_{1})$. Then either $8y=f(12y)$ or $8y=f(24y)$. But $f(12y)=4y$, since $4y \not\in B$, by the definition of $f$ on $B_{1}$. Otherwise $24y \in B$, in which case $9(9y)=3(24y)+(9y)$, 
provoking another contradiction.  
\par Thus we may define 
\begin{eqnarray*}
f(x)=8y, \;\;\;\; {\hbox{if $12y \in B$ and $3y \in B$}}.
\end{eqnarray*}
This completes the definition of $f$ on $B_{2}$, and it is automatic that $f(B_{1}) \cap f(B_{2}) = \phi$. 
\\
\\
{\em Step 3} : We show that $f \mid_{B_{2}}$ is one-to-one.
\\
\\
So suppose that there are $y,z \in A$ with $y \neq z$ but $f(9y) = f(9z)$. 
Without loss of generality, 
$f(9y)/9y < f(9z)/9z$. We then have nine cases to consider.
\\
\\
{\sc Case I} : $f(9y)=y$, $f(9z)=2z$ and $12y \not\in B$. 
\\
\\
Then $y=2z$ and $\{2y,6y\} \subseteq B$. But $2y=4z$ and then the equation $9(4z)=3(9z)+(9z)$ contradicts $B$'s avoidance of ${\cal L}_{3}$.
\\
\\
{\sc Case II} : $f(9y)=y$, $f(9z)=2z$ and $12y \in B$. 
\\
\\
Then $12y=24z$ and the equation $9(9z)=3(24z)+(9z)$ contradicts $B$'s avoidance
of ${\cal L}_{3}$. 
\\
\\
{\sc Case III} : $f(9y)=y$, $f(9z)=4z$ and $12y \not\in B$.
\\
\\
Then $6y=24z \in B$ and we get the same contradiction as in Case II.
\\
\\
{\sc Case IV} : $f(9y)=y$, $f(9z)=4z$ and $12y \in B$.
\\
\\
Then we still have that $6y=24z \in B$, so we get the same contradiction as in Case II.
\\
\\
{\sc Case V} : $f(9y)=y$, $f(9z)=8z$ and $12y \not\in B$.
\\
\\
Then $12z = 3y/2 \in B$ and the equation 
$9(\frac{3y}{2}) = 3(\frac{3y}{2})+(9y)$ yields a contradiction.
\\
\\
{\sc Case VI} : $f(9y)=y$, $f(9z)=8z$ and $12y \in B$.
\\
\\
Then we still have that $12z = 3y/2 \in B$, so we get the same contradiction as in Case V.
\\
\\
{\sc Case VII} : $f(9y)=2y$ and $f(9z)=4z$.
\\
\\
Then $y=2z$ and $12y=24z \in B$, so we get the same contradiction as in Case II.
\\
\\
{\sc Case VIII} : $f(9y)=2y$ and $f(9z)=8z$.
\\
\\
Then $12z = 3y \in B$, contradicting the definition of $f$ and the fact that $f(9y)=2y$. 
\\
\\
{\sc Case IX} : $f(9y)=4y$ and $f(9z)=8z$.
\\
\\
Then $y=2z$ and $3z=3y/2 \in B$, so we get the same contradiction as in Case V.
\\
\\
We have now completed Steps 1,2 and 3, and with that the proof of 
\\ Theorem 2.8.

\section{Results for Family II equations}

In this section, $\cal L$ denotes an equation $b(x+y) = cz$, where 
$b$ and $c$ are positive integers such that $b > 1$
and GCD$(b,c) = 1$. These equations were briefly 
touched on in [4], and the case $b=1$ was studied in detail in [1]. 
We present a theorem which describes extremal 
$\cal L$-avoiding subsets of $[1,n]$ for all values of $b,c$ and $n$. The
most interesting part of the theorem is part (i) which shows that the 
situation when $c >> b$ is quite different from when $b=1$, since the
extremal sets we describe are a $\lq$hybrid' between the two 
possibilities predicted by the Question in the introduction.  
\\
\\
{\bf Theorem 3.1}
{\em (i) If $c > 2b$ then, for every $n > 0$, the set
\begin{eqnarray*}
A_{n} = \left( {2bn \over c}, n \right] \cup 
\left\{ x \in \left[1, {2bn \over c}
\right] : {\hbox{$x \not\equiv 0 \; ({\hbox{mod $b$}})$}} \right\}
\end{eqnarray*}
is an $\cal L$-avoiding subset of $[1,n]$ of maximum size. 
\\
(ii) If $2 \leq c < 2b$ then, for every $n > 0$, the set 
\begin{eqnarray*}
A^{\prime}_{n} = \left\{ x \in [1,n] : x \not\equiv 0 \; ({\hbox{mod $b$}}) 
\right\}
\end{eqnarray*}
is an $\cal L$-avoiding subset of $[1,n]$ of maximum size.
\\
(iii) If $a = 1$ then, for every $n > 0$, the set 
\begin{eqnarray*}
A^{\prime \prime}_{n} = \left( {n \over 2b}, n \right]
\end{eqnarray*}
is an $\cal L$-avoiding subset of $[1,n]$ of maximum size.}
\\
\\
Note that 
\begin{eqnarray}
|A_{n}| = n - \lfloor {2n \over c} \rfloor, \\
|A^{\prime}_{n}| = n - \lfloor {n \over b} \rfloor, \\
|A^{\prime \prime}_{n}| = n - \lfloor {n \over 2b} \rfloor.
\end{eqnarray}
{\sc Proof of part (i)} : 
We fix $c > 2b$ and proceed by induction on $n$. The theorem
obviously holds if $n=1$. Fix $n > 1$ and let
$B$ be any $\cal L$-avoiding subset of $[1,n]$. We must show that 
$|B| \leq |A_{n}|$. 
\\
\\
First suppose there exists a number $z \in B \cap 
(\frac{bn}{c},\frac{2bn}{c}]$ which is a multiple of $b$. 
Let $z = bz_{1}$. Then $t := cz_{1} \in 
(n,2n]$ and, since $B$ avoids $\cal L$, there are no solutions in $B$ to 
the equation 
\begin{eqnarray*}
x+y=t.
\end{eqnarray*}
Now the map $f : x \mapsto t-x$ is a 1-1 mapping from the interval 
$I := [t-n,n]$ to itself, and for each $x \in I$, at most one of $x$ and 
$f(x)$ lies in $B$. Define $s \in [1,n]$ so that $t-n = n-s+1$. Then we
conclude that 
\begin{eqnarray*}
|B \cap [n-s+1,n]| \leq \lfloor s/2 \rfloor.
\end{eqnarray*}
If $s=n$ then $|B| \leq \lfloor \frac{n}{2} \rfloor \leq |A_{n}|$, by (30). 
Otherwise, the induction hypothesis yields that 
\begin{eqnarray*}
|B| \leq |A_{n-s}| + \lfloor {s \over 2} \rfloor \\ = 
(n-s) - \lfloor {2(n-s) \over a} \rfloor + \lfloor {s \over 2} \rfloor \\
< |A_{n}| + 1, 
\end{eqnarray*}
and since $|B|$ is an integer, we conclude that $|B| \leq |A_{n}|$, as desired.
\\
\\
Thus we may assume that $B$ contains no multiples of $b$ in the interval 
$(\frac{bn}{c},\frac{2bn}{c}]$. If $B$ contains no multiples of $b$ at all in
the range $[1,\frac{2bn}{c}]$, then trivially $|B| \leq |A_{n}|$. So let's 
assume $B$ does contain such a number, and let the largest such be 
$z_{0}$. Thus $z_{0} \in [1,\frac{bn}{c}]$. Let $z_{0} = bz_{1}$. Then 
$cx_{1} \in [c,n]$ and there are no solutions in $B$ to 
\begin{eqnarray}
x+y = c z_{1}.
\end{eqnarray}
Note that since $c > 2b$, we have $cz_{1} > 2z_{0}$, so if $(x,y)$ is a 
solution to (33), then $x \leq z_{0} \Rightarrow y > z_{0}$. 
\par Suppose $cz_{1} \equiv j \; ({\hbox{mod $b$}})$. Then for every number
$x \in [1,z_{0}]$ such that $x \not\equiv j \; ({\hbox{mod $b$}})$, at most 
one of $x$ and 
$cz_{1} - x$ lies in $B$. But note 
that, for every such $x$, $cz_{1} - x$ is not divisible by $b$, and is 
strictly greater than $z_{0}$. We conclude that at least 
$z_{0} - z_{1} \geq z_{1}$ numbers are missing from $B$ which are either 
multiples of $b$ in $[1,z_{0}]$ or not divisible by $b$. 
This implies that $|B| \leq |A_{n}|$ and 
completes the proof of part (i) of Theorem 3.1.
\\
\\
{\sc Proof of Part (ii)} : We divide the proof into two cases.
\\
\\
{\em Case I} : $2 \leq c < b$.
\\
\par Fix $n > 0$ and an $\cal L$-avoiding 
subset $B$ of $[1,n]$. We must show that 
$|B| \leq |A^{\prime}_{n}|$. If $B$ contains no multiples of $b$ 
we are done, so suppose the contrary. Let $z = bz_{1}$ be the largest element 
of $B$ which is 
a multiple of $b$. Then it suffices to produce at least $z_{1}$ numbers in the
interval $[1,z]$ which are not in $B$. Since $B$ avoids $\cal L$, 
it contains no 
solutions to the equation 
\begin{eqnarray*}
x+y = cz_{1}.
\end{eqnarray*}
Thus $B$ contains no more than $\lceil \frac{cz_{1}}{2} \rceil$ of
the numbers in the interval $[1,cz_{1}]$. 
But since $2 \leq c < b$, it follows that $z_{1} \leq \lfloor 
\frac{cz_{1}}{2} \rfloor$ and 
$cz_{1} < z$. Thus we are done.
\\
\\
{\em Case II} : $b < c < 2b$.
\\
\par We proceed by induction on $n$. 
The theorem obviously holds for $n=1$. Now 
fix $n > 1$ and an $\cal L$-avoiding 
subset $B$ of $[1,n]$. If $B$ contains no multiples 
of $b$ then we are done, so we may assume that $B$ contains 
some such elements. 
\par First suppose there exists $z \in B \cap (\frac{b}{c} n,n]$ 
which is a multiple of $b$. Let $z = bz_{1}$. 
Then $cz_{1} \in (n,\frac{c}{b}n] \subseteq (n,2n]$. 
To simplify notation, set $cz_{1} := t_{1}$ and $t_{1}-n := n-s_{1}+1$ . 
Since $B$ avoids $\cal L$, it contains no solutions to the equation 
\begin{eqnarray*}
x+y=t_{1}.
\end{eqnarray*}
The map $f : x \mapsto t_{1}-x$ is a 1-1 mapping from the interval 
$I_{1} = [n-s_{1}+1,n]$ to itself and for each $x \in I_{1}$ at most one of 
the numbers $x$ and $f(x)$ lies in $B$. Thus 
\begin{eqnarray*}
\mid B \cap I_{1} \mid \; \leq \lfloor s/2 \rfloor.
\end{eqnarray*}
Then clearly (since $b \geq 2$) the induction argument implies 
that $|B| \leq |A^{\prime}_{n}|$. 
\par Thus we may assume that $B \cap (\frac{b}{c} n,n]$ contains no 
multiples of 
$b$. But then another application of the induction hypothesis yields 
that $|B| \leq |A^{\prime}_{n}|$ in this case too. 
\par Thus part (ii) of the theorem is proved.    
\\
\\
{\sc Proof of Part (iii)} : The argument is similar to that in Case I of part 
(ii). Let $B$ be an $\cal L$-avoiding 
subset of $[1,n]$. If $B$ contains no multiples of
$b$, then clearly $|B| \leq |A^{\prime \prime}_{n}|$. So suppose 
$z = bz_{1}$ is the largest multiple of $b$ in $B$. Then 
\begin{equation}
|B \cap (z,n]| \leq (n-z) - \lfloor {n-z \over b} \rfloor.
\end{equation}
Since $B$ avoids $\cal L$, it contains no solutions to the equation 
\begin{eqnarray*}
x+y = z_{1}.
\end{eqnarray*}
Thus
\begin{equation}
|B \cap [1,z]| \leq z - \lfloor {z_{1} \over 2} \rfloor = z - \lfloor
{z \over 2b} \rfloor.
\end{equation}
Clearly, (34) and (35) imply that $|B| \leq |A^{\prime \prime}_{n}|$ (with 
strict inequality unless $z_{1} = \lfloor n/b \rfloor$).
\par This completes the proof of Theorem 3.1.
\\
\\
{\bf Concluding remark} It is worthwhile to investigate if the proof
of Theorem 3.1 can be used to obtain a stronger result, namely a 
classification of the extremal sets. We choose not to go into this matter in 
this paper, which we think already contains enough in the way of detailed, 
technical computations. In any case, the important thing is the 
$\lq$hybrid' nature of the extremal sets in part (i) of the theorem. Note that,
for each $b \geq 2$, 
the sets $A_{n}$ have strictly greater asymptotic 
density than any sets of type {\bf II}$^{\prime}$. This follows from (30),
Lemma 1(b) of [4] and a straightforward computation. 
    
\section*{References}

[1] {\sc A. Baltz, P. Hegarty, J. Knape, U. Larsson and T. Schoen}, The 
structure of maximum subsets of $\{1,...,n\}$ with no solutions to $a+b=kc$,
{\em Electron. J. Combin.} {\bf 12} (2005), Paper No.19, 16pp.
\newline
[2] {\sc F.R.K. Chung and J.L. Goldwasser}, Integer sets containing no 
solutions to $x+y=3z$, in : R.L. Graham and J. Ne$\check{s}$et$\check{r}$il 
eds.,
The Mathematics of Paul Erd\H{o}s, Springer, Berlin (1997), pp. 218-227.
\newline
[3] {\sc J.M. Deshouillers, G.A. Freiman, V. S\'{o}s and M. Temkin}, On
the structure of sum-free sets 2, in : J.M. Deshouillers et al (eds.), 
Structure Theory of Set Addition, {\em Ast\'{e}risque} {\bf 258} (1999), 
149-161.
\newline
[4] {\sc K. Dilcher and L.G. Lucht}, Finite pattern-free sets of integers, 
{\em Acta Arith.} {\bf 121}, No.4, (2006), 313-325.
\newline
[5] {\sc J. Knape and U. Larsson}, Sets of integers and permutations avoiding 
solutions to linear equations, Master's Thesis, G\"{o}teborg
University, 2004. Available online at
http://www.mdstud.chalmers.se/$\sim$md0larur/sista$\underline{\;}$magex.ps
\newline
[6] {\sc L.G. Lucht}, Dichteschranken f\"{u}r die L\"{o}sbarkeit gewisser
linearer Gleichungen, {\em J. Reine Angew. Math.} {\bf 285} (1976), 209-217.
\newline
[7] {\sc L.G. Lucht}, Extremal pattern-free sets of positive integers,
{\em Ann. Univ. Sci. Budapest, Sect. Comp.} {\bf 22} (2003), 253-268.
\newline
[8] {\sc T. \L uczak and T. Schoen}, On infinite sum-free sets of natural
numbers, {\em J. Number Theory} {\bf 66} (1997), 211-224.
\newline
[9] {\sc T. \L uczak and T. Schoen}, Solution-free sets for linear 
equations, {\em J. Number Theory} {\bf 102} (2003), 11-22. 
\newline
[10] {\sc I.Z. Ruzsa}, Solving a linear equation in a set of integers I, 
{\em Acta Arith.} {\bf 65} (1993), 259-282.
\newline
[11] {\sc I.Z. Ruzsa}, Solving a linear equation in a set of integers II, 
{\em Acta Arith.} {\bf 72} (1995), 385-397.
\newline
[12] {\sc T. Schoen}, On sets of natural numbers without solution to a 
non-invariant linear equation, {\em Acta Arith.} {\bf 93} (2000), 149-155.

\end{document}